\documentclass[11pt, oneside]{article}  
\usepackage{geometry}                	
\usepackage{amssymb}
\usepackage{graphicx}
\usepackage{xcolor}
\usepackage{hyperref}

\def\int#1{{\mbox{\em int}}(#1)} 

\def\vp{\vspace{3ex}}

\usepackage[font={small}]{caption}

\title{ Notes on a strongly aperiodic monotile in $E^3$}
\makeatletter
\author{\small Chaim Goodman-Strauss\\\small  National Museum of Mathematics\\\small  chaimgoodmanstrauss@gmail.com}
\date{}
\begin{document}
\maketitle

\abstract{\small We provide a clearer presentation of the Chair44 monotile, a new strongly aperiodic monotile in three-dimensional Euclidean space $E^3$.}

\vp
\vp


 With LLM tools (ChatGPT Astra),
Ioannis Tsiokos recently reported the first known strongly aperiodic monotile in three-dimensional Euclidean space $E^3$, the {\em Chair44 monotile}. (The 44 is of archetypically marginal relevance. See Figure~\ref{fig4}, taken from Tsiokos' preprint~\cite{ioannis}. We give definitions in Section~\ref{sectbackground}.)
 Hats off to Tsiokos for this discovery, the framework he developed, and the human inquiry that led to the Chair44 monotile. 

\vp
However, the paper itself and the formalization that it claims for support are  not in a form that is readily usable to anyone who wants to understand or check it. This is  not just because of the needless length or complexity of the argument (which is indeed overwhelming as presented), or because the formal repository  is not in normal form, or because much of the text is misleading or poorly organized, but mainly because there is so much irrelevant puff to push through. 
We must, absolutely, insist on a higher standard for scientific discourse. The aim must be to communicate, clearly, with people in the community --- certainly this is an historical requirement for publication.

Even more recent events have only made the importance of this clearer. Problem \#155 in the recent Open AI release~\cite{openai_monotile} provides a counterexample to the Periodic Tiling Conjecture in three dimensions. This is  a weakly aperiodic translation monotile  derived from  Greenfeld and Tao's construction in much higher dimensions~\cite{greenfeld_tao}. (Greenberg and Tao's techniques use the product of a large number of  cyclic spaces; this new paper seems to encode this into a single large cyclic space using Sunzi's Remainder Theorem.) 
This presentation is in more proper form, but  its notation slides around and its big-picture descriptions aren't helpful. Its deficiencies give us a sense of the abilities of these machines in the fall of 2026  when used  with substantial computational resources and expertise for a prestige project but in bulk with little human attention.

No doubt this technology will rapidly improve, and any assessment we could make here is soon going to be ancient history. But our own meaning and understanding, for our own curiosity and taste, must always come first.\footnote{
Terry Tao has prepared a statement along these lines, signed by twenty-five Fields Medalists~\cite{tao}, and Tao's blog has become a rallying point for many writers. See too
  the Leiden Declaration, endorsed by the International Mathematical Union~\cite{leiden}.}

\vp Indeed, taking the interest to simplify the construction, we clarify it: Once redrawn (Figure~\ref{newl}), this monotile turns out to be remarkably elegant,  with a reasonable claim for the simplest-known ``Berger-style'' proof of aperiodicity~\cite{berger} (right, Figure~\ref{newl}).  

Felix Flicker has shed further light on the Chair44 (including a key observation that the markings we show here as blue arrows can be removed without changing the combinatorics of the tile, as well as an extensive bibliography and further  context)~\cite{flicker}. We find Pieter Mostert's 3d printable design to be very useful~\cite{mostert} (Figure~\ref{mosterttile}) and you can design your own~\cite{gs_ltile}. Within our tiling community, we have exchanged hundreds of messages exploring this tile and its variations, and using LLMs to do so. We can expect that Tsiokos' question and the LLM's discovery will lead to even more insight and investigation.  

If we argue and push for it,
this 
 may point in a positive direction: New tools in human hands together  with human insight produce new mathematics of   interest to humans.

\section{Background}\label{sectbackground}
Copies of the Chair44 monotile can be fitted together without gaps or overlaps to cover all of space (the shape is a {\em tile} that can form {\em tilings}, or {\em tessellations}, and is therefore a {\em monotile}). We allow rotations but do not permit reflected copies. (A {\em translation} monotile would allow only translated copies.) The Chair44 has the special property that any tiling that it forms  must be {\em non-periodic}, with no translational symmetry (and as it can occur in only finitely many orientations,  it therefore has no infinite-cyclic symmetry of any kind). The Chair44 monotile is {\em strongly aperiodic.}

(In contrast,
a tile or set of tiles is {\em weakly aperiodic} if it admits tessellations, but no periodic tessellation has a compact fundamental domain. The latter term was coined, in part, to describe the Schmitt-Conway-Danzer biprism, a monotile in space that doesn't allow tilings with a compact fundamental domain, but does allow tilings with an infinite screw symmetry~\cite{danzer}; similarly the recent counter-examples to the Periodic Tiling conjecture are weakly but not strongly aperiodic. In the plane, there is no need to distinguish between {weak} and {strong} aperiodicity: 
See Theorem 3.7.1 in~\cite{grsh}.)

\begin{figure}[htbp]
\centerline{
\includegraphics{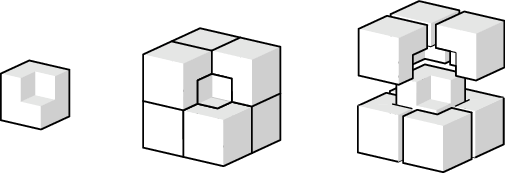}}
\caption{A three dimensional L-tile and a supertile; these generalize to all $E^{n\geq 2}$~\cite{gs_en}.}\label{L3}
\end{figure}

\begin{figure}[h]
\centerline{
\includegraphics{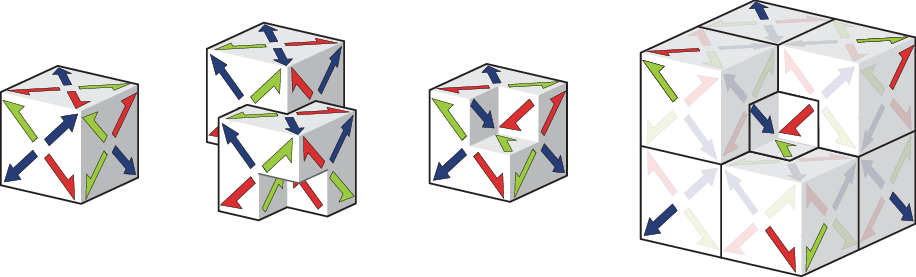}}
\caption{Several views of a re-rendering of  the Chair44 monotile~\cite{ioannis}. (Compare to Figure~3.)
 Blue arrows must match; green arrows must match with red ones; tiles fit together as second to left. In any tessellation satisfying these rules each tile must lie in a copy of the supertile shown at right.  In turn these supertiles fit together in just the same way that the original tiles can, completing a ``Berger-style" proof.}
\label{newl}
\end{figure}

\begin{figure}[h]
\centerline{
\includegraphics[width=.7\textwidth]{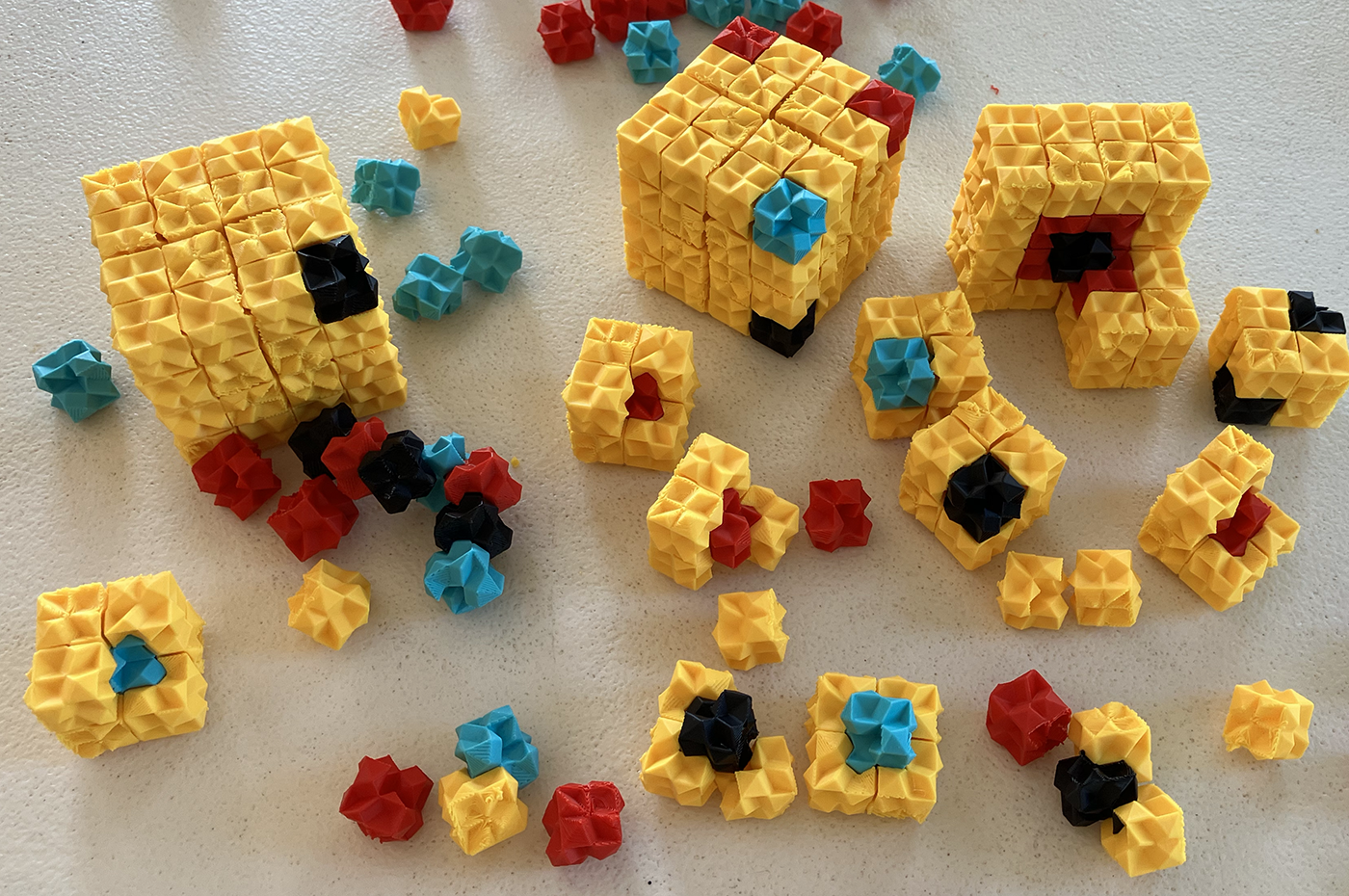}}
\caption{A pile of printed Chair44 monotiles, in a version designed by Pieter Mostert. You can download~\cite{mostert} or design~\cite{gs_ltile} your own!}
\label{mosterttile}
\end{figure}

  The construction is based on  the {\em three-dimensional L-tile}, shown at left in Figure~\ref{L3},  made by deleting a half-sized cube from the corner of a large one.   
  This shape is a {\em rep-tile} which means copies of it can be fitted together to form a larger copy of itself. 
  We show this larger cluster of tiles, which we call a {\em supertile}, in the middle of the same figure, and at right in exploded form. In the same way that tiles may  form supertiles, supertiles can be fitted together  into a 2-level supertile, and 2-level supertiles into 3-level ones, and so on {\em ad infinitum}. 

It's  believable and true, but more subtle than it seems (see Theorem 3.8.1 of Gr\"unbaum and Shephard's masterwork~\cite{grsh}) that if we can form arbitrarily high-level supertiles then we can  form infinite tilings, tilings that seem to be infinite-level supertiles. 
These tilings will be {\em hierarchical}:  each tile is within a supertile (a unique one), which is within a (unique) 2-level supertile, and so on. Each tile is within an infinite hierarchy of supertiles, and this hierarchy is unique. Such a hierarchical tiling must be {\em non-periodic}, and cannot have any translational symmetry at all.  On the other hand, this $L$-tile shape can also form periodic tessellations or more random structures. It is not aperiodic. 

Over sixty years ago, Robert Berger realized hierarchical tilings like this one could be the basis for constructing {aperiodic}  tiles~\cite{berger}.  The trick is to design restrictions (like the bumps and indentations on jigsaw puzzle pieces, or marked colors that must match) that allow the tiles to fit together into these hierarchical structures, but {\em only} allow that.  

\vp
 Today there are dozens of known examples of sets of tiles that ``enforce" rep-tile like rules in this way, mostly found by some sort of magic. However, every known example has been more complex  --- either the tiles needed to be marked in more than one way to ensure they could only form these hierarchical structures, or the  rep-tile rules were more complex.\footnote{The spectre and the tiles in the hat tile continuum have rep-tile like rules,~\cite[Figure 2.11]{hat} and~\cite[Figure 2.1]{spectre}, but strangely each has {\em two} rules on two different shapes, fractals that bear no outward resemblance to the original tiles.}
 No one had yet found an aperiodic monotile enforcing a true rep-tile rule, or one in three-dimensions.  For almost up-to-date histories, we refer the reader to the papers \cite{hat} and~\cite{spectre} announcing two-dimensional aperiodic monotiles, as well as Flicker's preprint on the Chair44, with its extensive bibliography~\cite{flicker}.
 
 \vp
 Tsiokos discovered the Chair44 monotile through  prompts to an LLM, together with his  {\em Six-birds} framework, further referenced in~\cite{ioannis}.
 Tsiokos and the LLM give an account of the process in their Section~3.3. The LLM proposed using the three-dimensional L-tile, as it is a rep-tile in space, and thus is a good candidate for an aperiodic monotile.\footnote{
The shape and its substitution rule generalize to all dimensions. 
In previous work~\cite{gs_en},  we see that for each $n\ge 3$, there is an aperiodic {\em pair} of tiles that enforces this rep-tile tiling in $n$-dimensional space. One of those tiles has disconnected interior but it seems likely that tile can be reconfigured to have connected and simply connected interior.}

After a false start, they employed a  sensible search procedure: They enumerated all the ways that neighboring L-tiles meet within supertiles, 2-level supertiles and so forth, and consequently the ways that neighbors meet within a full hierarchical L-tiling. From this (as I understand), they constructed a maximal set of matching rules, restrictions placed upon the ways the tiles might fit together, that still would allow a hierarchical tiling.

They then exhaustively checked the ways that copies of this marked tile could fit together to see that neighborhoods of this marked tile, the Chair44 monotile, could only resolve into pieces of tilings by supertiles --- showing that any tiling by the marked tiles could  be grouped into a tiling by the supertiles, and that the supertiles could only fit together in the same way the tiles do. By 
Berger's argument, which we demonstrate in Section~\ref{sect2} and illustrate in Figure~\ref{newl},  the tile is aperiodic. 

There's no {\em a priori} reason that this will pan out, that the rules that are loose enough to allow the configurations that appear in the hierarchical L-tilings are also restrictive enough to allow nothing else. 
But if there's something there, this method will find it. In that light, it is natural that this process found an especially  elegant aperiodic monotile:  if it exists to be found,  it should be among the first to be discovered in this way. 
For us in the tiling community, this is a neat result. Even better, this points to the use of new LLM tools for many more discoveries  along these lines. We can hope for more!\footnote{Joseph Myers suggests that aperiodic monotiles might naturally be  easier to discover in higher dimensions.} 

\begin{figure}[h]
\centerline{
\includegraphics[scale=.75]{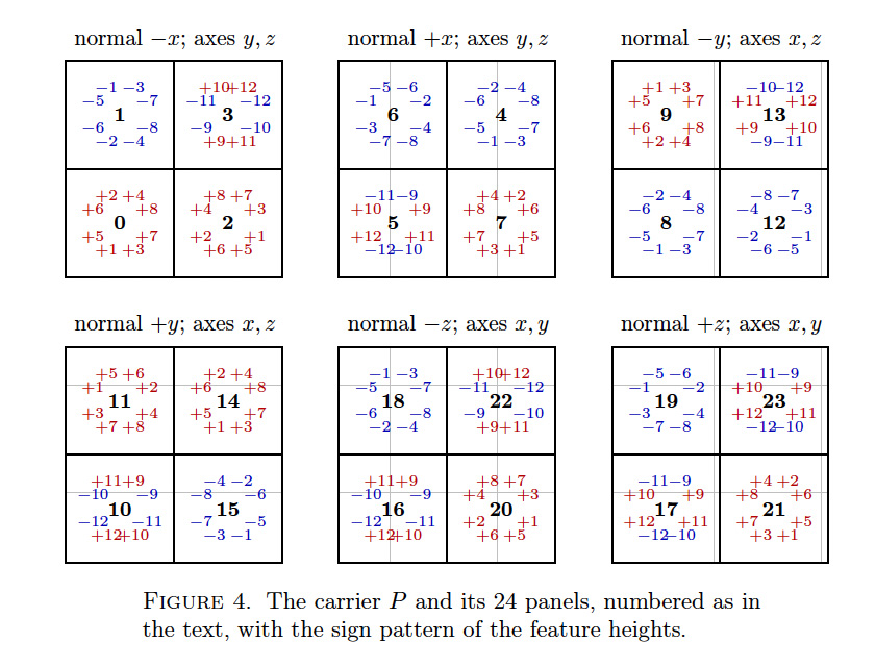}}
\caption{Figure 4 of Tsiokos's preprint announcing the Chair44 monotile~\cite{ioannis}.}
\label{fig4}
\end{figure}

However, we must  insist on more transparent and helpful presentation of the results:
Once we recognize that the Chair44 monotile is a marked three-dimensional L-tile, and which way it is pointed, and that it will be used as a rep-tile, the only  information useful for understanding it, out of all of~\cite{ioannis}, appears in that paper's Figure 4, reproduced here as Figure~\ref{fig4}. There, in some form,  we find its markings. (Compare to our Figure~\ref{newl}.) 

These  markings are shown from the six coordinate directions, but in half (which three?) the view is reversed from how it would appear to someone outside of the tile. It takes some effort to work out where the missing cube is supposed to be.   The  twenty-four strangely-ordered panel numbers  aren't  useful to a human reader, and at a glance we  expect there are much simpler markings (the small numbers are the heights of little pyramids).  

\vp 
As it turns out,  the core idea is very simple, and a much shorter paper is appropriate: Any reader can cut out and assemble the plans in Figure~\ref{net}, to have an intuitive sense of how copies of the Chair44 monotile fit together.  Some readers will be able to visualize the structure well enough in their minds to verify the claims below. Others may wish to draw or model the tile for themselves.
In my experience,  people grasp the core combinatorial argument within a few minutes of handling printed plastic copies of the monotile~\cite{mostert}.  As people, we can verify such claims with our mental toolkit --- we offload the computational complexity of the proof into our hands, eyes and minds.

\section{The tile}  \label{sect2}

We've already met the three-dimensional L-tile in Figure~\ref{L3}, and its rep-tile rule, pulled apart at right in that figure so we can better see how the tiles fit together: We can see what we'll call a ``central'' L-tile in the middle, and seven ``outer'' tiles surrounding it. 

\begin{figure}[h]
\centerline{
\includegraphics[scale=.9]{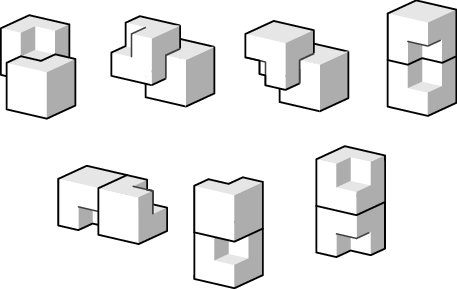}}
\caption{Adjancies between L-tiles in $n$-level supertiles~\cite{gs_en}.}
\label{adj}
\end{figure}

\begin{figure}[h]
\centerline{
\includegraphics[width=.9\textwidth]{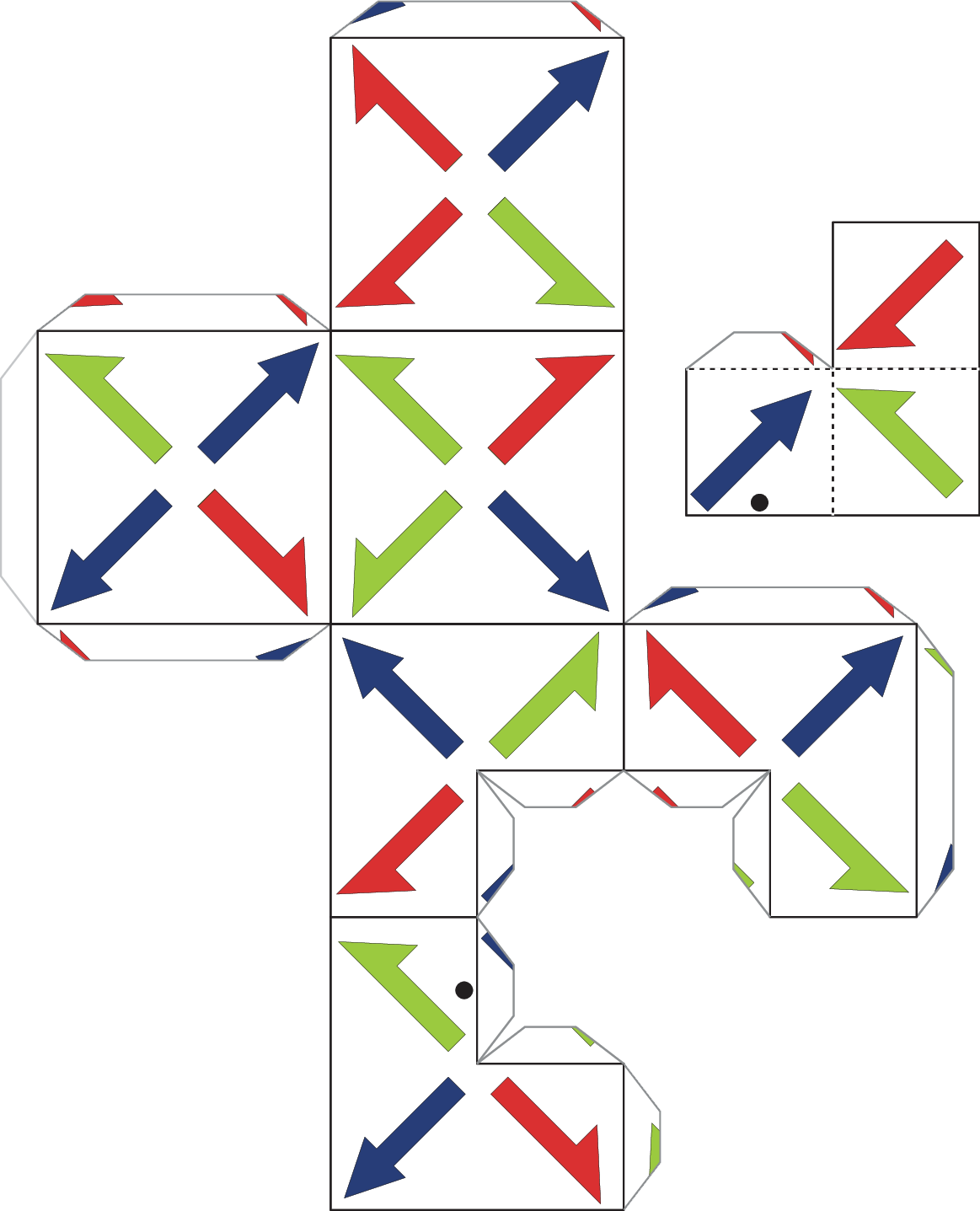}}
\caption{Cut out and assemble to make your own! (All fold lines are ``mountain" folds, except the two dashed lines which are `valley" folds. The black dots will be adjacent across an edge of the model. The tabs can go on the inside or the outside of the model.)}
\label{net}
\end{figure}

 In Figure~\ref{adj}, taken from the paper~\cite{gs_en} which generalizes this to all higher dimensions, we show the ways that L-tiles may meet along a face when the rep-tile rule is iterated, corresponding with the 44 positions for which the new tile is named~\cite{ioannis}: At top left we show how two outer tiles in a supertile meet; along the rest of the top we show how the inner tile in a supertile may meet an outer one. When we fit  supertiles together in one of those three ways, we'll find a pair of tiles meeting the same way as left on the bottom row (shown in a different orientation at middle). When supertiles are placed into any other of these arrangements, we find tiles arranged as at bottom right. 
These are all the possible ways tiles could be adjacent in a hierarchical tiling --- no new arrangements will come about in putting supertiles together in ways tiles can fit together. 

In Figure~\ref{newl} we rework the markings from Figure~\ref{fig4}, showing the tile from three angles. At right, eight of these tiles have been assembled into a marked supertile, with some markings lightened to better see the structure.

There are three kinds of markings, shown as a coloring: red, green and blue\footnote{Flicker notes that the ``blue'' markings are redundant and can be removed entirely~\cite{flicker}.}(which will appear medium gray, light gray and dark gray if printed in black and white). Red may only meet green, and green may only meet red; blue may only meet blue. Notice that at each of the seven ``corners'' (what we'll call the  vertices of the original cube that remain on the L-tile) as well as at the ``front socket''  (the concave vertex at the center of the original cube), we find each of the three markings, always in the same cyclic ordering. The active reader is invited to make the correspondence with the markings in Figure~\ref{fig4} explicit.  

To disambiguate the orientation of the markings on fronts and back, the tile fits together with a translation of itself, its front socket fitting into its ``back'' corner (the corner opposite the corner that was removed from the cube to form the L-tile).\footnote{For a moment it may be confusing that the same cyclic ordering of the colors can match at the front socket and an outer corner --- but when we view the scene from one direction, we  flip one tile over, which reverses the order of the markings from our point of view, and green and red must match, which reverses the order again.}  We show this second from left in Figure~\ref{newl}.

By inspection, a tile can never fit with its mirror reflection, so it does not matter whether or not we allow the tiles to be reflected (which we cannot do with physical three-dimensional tiles in any case).  The tiles are chiral, and its tilings (once we show there are any) are homochiral~\cite{spectre}: within a tiling  all copies of the tile have the same orientation.   

Consequently we can create a version of the tile that is an unmarked shape, converting the markings into bumps and nicks. The blue marking can be raised on one side of a diagonal and correspondingly lowered on the other (or, as Flicker notes, simply left flat). Designs for the red and green markings need only to fit into each other. 

With physical models in hand, the active reader needs only a few moments to  verify that any marked tile is either at the corner or the center of a marked supertile configuration like that of Figure~\ref{L3}. \footnote{We have given 3d printed copies of Mostert's Chair44 monotile to people with no instructions. Very quickly on their own, they seem to  discover the configuration at right in Figure~\ref{newl}, and to realize that's the only way to pack them together without gaps, completing the combinatorial core of this proof. The rest is packaging. In a sense this calculation is just as complex as a formal proof that checks how cubes fit together (a component of~\cite{ioannis}), but we are able to make use of our embodied physical sense of space. You can print your own copies of Mostert's Chair44~\cite{mostert} or one of your own design~\cite{gs_ltile}.} (Note that the inside socket of any tile must meet the outer corner of some tile. If this corner is not the back corner, then the first tile is one of seven surrounding the second, which is at the center of a supertile. If the corner is the back corner, one or the other of the first or the second tiles will be at the center of a supertile.)

The corners of this supertile are each the back corner of one of the outer L-tiles that form it. Each back corner of a tile matches its inner socket, with red and green swapped, which matches the corner of the inner L-tile it meets, with red and green swapped a second time. 
Therefore,  the markings on the  corners of the supertile are the same as the markings on the  corners of the inner tile! 

At its corners a supertile's markings are the same as  those at the corners of a marked tile.  
The rest of the markings do not vary from face to face (see at right in Figure~\ref{newl}). Supertiles can only meet fully face-to-face, and they may only meet in the same way that the tiles do. They can form marked 2-level supertiles. Moreover, in any tiling by the Chair44 monotile, they must. 

This continues. At each level, the supertiles are geometrically just larger three-dimensional $L$-tiles, and as shapes can only fit together in the way that those do, and the markings at their corners can only match in the way that the markings of the Chair44 monotile can. Each level of marked supertile can be assembled into a next-level marked supertile, and in any tiling by the Chair44 monotile, each tile must lie in a unique hierarchy of supertiles up to that next level. Repeating this {\em ad infinitum}, any tiling by the Chair44 monotile is hierarchical, and cannot admit any translational symmetry. We've completed a  Berger-style proof of:

\vp{\noindent\bf Theorem:} The tile of Figure~\ref{newl}, a redrawing of the Chair44 monotile~\cite{ioannis}, is a strongly aperiodic monotile.

\section{Acknowlegements} Thanks to  Craig S. Kaplan,  Felix Flicker, Pieter Mostert,  V. Van Dongen, and members of the Tiling Group for helpful comments and insights; and to  Shosha Wheeler, 
Max Grossman, 
Macy Aiken, 
Susanne Goldstein, 
Anthony Wright, 
Irene Yang, and
Jade Nichols  for play-testing Figure~\ref{net}. 
No LLMs were used in the preparation of this paper.


\begin{thebibliography}{6}


\bibitem{berger} {R. Berger}, {\em The undecidability of the Domino Problem},
Memoirs Am. Math. Soc.  {\bf 66} (1966).



\bibitem{danzer} {L. Danzer}, {\em A family of 3D-spacefillers not permitting
any periodic or quasiperiodic tiling}, Aperiodic `94, World 
Scientific (1995), 11-17.

\bibitem{flicker} F. Flicker, {\em Matching rules for a three-dimensional strongly aperiodic monotile}, \url{arxiv.org/abs/2609.23783}.

\bibitem{gs_en} {C. Goodman-Strauss}, {\em An aperiodic pair of 
tiles in  $E^n$
for all
$n\geq 3$},  Europ. J. Combinatorics {\bf 20} (1999), 385-395.


\bibitem{gs_ltile}{C. Goodman-Strauss}, {\em Chair44 viewer}, \url{chaimgoodmanstrauss.com/ltile}









\bibitem{greenfeld_tao} R. Greenfeld and T. Tao, {\em A counterexample to the periodic tiling conjecture},  Annals of Mathematics {\bf  200} (2024), pages 301–363.

\bibitem{grsh} {B. Gr\"unbaum} and {G.C. Shepherd}, {\em Tilings and Patterns}, 2nd ed., Dover (2016). 

\bibitem{leiden} ``The Leiden Declaration,"
{\em leidendeclaration.ai}

\bibitem{mostert} P. Mostert, placed at \url{chaimgoodmanstrauss.com/downloads/chair44_mostert.stl}




\bibitem{openai_monotile} OpenAI, {\em A translational tile with no fully periodic tiling
in dimension three}, preprint at
\url{ github.com/openai/math/tree/main/preprints/A-translational-tile-with-no-fully-periodic-tiling-in-dimension-three-September-23-2026}


\bibitem{hat}
David Smith, Joseph Samuel Myers, Craig S. Kaplan, and Chaim Goodman-Strauss. 
``An aperiodic monotile.'' \textit{Combinatorial Theory}, vol. 4, no. 1, 2024.


\bibitem{spectre}
David Smith, Joseph Samuel Myers, Craig S. Kaplan, and Chaim Goodman-Strauss. 
``A chiral aperiodic monotile.'' \textit{Combinatorial Theory}, vol. 4, no. 2, 2024.

\bibitem{tao} T. Tao
``A Severe Misalignment of AI in Mathematics" \url{terrytao.wordpress.com/2026/09/11/a-severe-misalignment-of-ai-in-mathematics}


\bibitem{ioannis} I. Tsiokos, {\em A strongly aperiodic monotile in three dimensions}, arXiv:2609.19214v1

\end{thebibliography}
\end{document}